\documentclass[11pt]{amsart}
\usepackage{amssymb}
\usepackage{mathrsfs}
\usepackage{amsfonts}
\usepackage[all]{xy}

\usepackage{pstricks, pst-tree}
\allowdisplaybreaks

\begin{document}
\theoremstyle{plain}
\newtheorem{thm}{Theorem}[section]
\newtheorem{theorem}[thm]{Theorem}
\newtheorem{lemma}[thm]{Lemma}
\newtheorem{corollary}[thm]{Corollary}
\newtheorem{proposition}[thm]{Proposition}
\newtheorem{conjecture}[thm]{Conjecture}
\newtheorem{obs}[thm]{}
\theoremstyle{definition}
\newtheorem{construction}[thm]{Construction}
\newtheorem{notations}[thm]{Notations}
\newtheorem{question}[thm]{Question}
\newtheorem{problem}[thm]{Problem}
\newtheorem{remark}[thm]{Remark}
\newtheorem{remarks}[thm]{Remarks}
\newtheorem{definition}[thm]{Definition}
\newtheorem{claim}[thm]{Claim}
\newtheorem{assumption}[thm]{Assumption}
\newtheorem{assumptions}[thm]{Assumptions}
\newtheorem{properties}[thm]{Properties}
\newtheorem{example}[thm]{Example}
\newtheorem{comments}[thm]{Comments}
\newtheorem{blank}[thm]{}
\newtheorem{defn-thm}[thm]{Definition-Theorem}

\newcommand{\sM}{{\mathcal M}}


\title[The n-point functions for intersection numbers]{The n-point functions for intersection numbers\\ on moduli spaces of curves}
        \author{Kefeng Liu}
        \address{Center of Mathematical Sciences, Zhejiang University, Hangzhou, Zhejiang 310027, China;
                Department of Mathematics,University of California at Los Angeles,
                Los Angeles, CA 90095-1555, USA}
        \email{liu@math.ucla.edu, liu@cms.zju.edu.cn}
        \author{Hao Xu}
        \address{Center of Mathematical Sciences, Zhejiang University, Hangzhou, Zhejiang 310027, China}
        \email{haoxu@cms.zju.edu.cn}

        \begin{abstract}
           Using the celebrated Witten-Kontsevich theorem, we prove a recursive formula of the $n$-point functions for
           intersection numbers on moduli spaces of curves.
           It has been used to prove the Faber intersection number conjecture and motivated
           us to find some conjectural vanishing identities for Gromov-Witten invariants. The latter has been proved
           recently by X. Liu and R. Pandharipande. We also give a combinatorial interpretation of $n$-point functions in terms of summation over binary trees.
        \end{abstract}
    \maketitle

\maketitle

\section{Introduction}

 Let $\overline{\sM}_{g,n}$ be the moduli space of stable
$n$-pointed genus $g$ complex algebraic curves and $\psi_i$ the
first Chern class of the line bundle corresponding to the cotangent
space of the universal curve at the $i$-th marked point. Let
$\mathbb E$ denote the Hodge bundle. The fiber of $\mathbb E$ is the
space of holomorphic one forms on the algebraic curve. Let us denote
the Chern classes by
$$
\lambda_k=c_k(\mathbb E),\quad 1\leq k\leq g.
$$
More background material about moduli spaces of curves can be found
in the paper \cite{Vak}.

We use Witten's notation
$$\langle\tau_{d_1}\cdots\tau_{d_n}\rangle_g:=\int_{\overline{\sM}_{g,n}}\psi_1^{d_1}\cdots\psi_n^{d_n}.$$

These intersection numbers are the correlation functions of two
dimensional topological quantum gravity. Motivated by an analogy
with matrix models, Witten \cite{Wi} made a remarkable conjecture
(originally proved by Kontsevich \cite{Ko}) that the generating
function
\begin{equation}\label{gen}
F(t_0, t_1, \ldots)= \sum_{g} \sum_{\bold n}
\langle\prod_{i=0}^\infty \tau_{i}^{n_i}\rangle_{g}
\prod_{i=0}^\infty \frac{t_i^{n_i} }{n_i!}
\end{equation}
 is a $\tau$-function for
the KdV hierarchy, which also
provides a recursive way to compute all these intersection numbers.
In particular, $U=\partial^2 F/\partial t_0^2$ satisfies the classical KdV equation
\begin{equation} \label{KdV}
\frac{\partial U}{\partial t_1}=U\frac{\partial U}{\partial t_0}+\frac{1}{12}\frac{\partial^3 U}{\partial t_0^3}.
\end{equation}

Witten's conjecture was reformulated by Dijkgraaf, Verlinde, and
Verlinde [DVV] in terms of the Virasoro algebra. Now there are
several new proofs of Witten's conjecture \cite{CLL, KaL, KiL, Mi,
OP}.

\begin{definition}
We call the following generating function
$$F(x_1,\dots,x_n)=\sum_{g=0}^{\infty}\sum_{\sum d_j=3g-3+n}\langle\tau_{d_1}\cdots\tau_{d_n}\rangle_g\prod_{j=1}^n x_j^{d_j}$$
the $n$-point function.
\end{definition}

The $n$-point function is an alternative way to encode all
information of intersection numbers of $\psi$ classes. Okounkov
\cite{Ok} obtained an analytic expression of the $n$-point functions
in terms of $n$-dimensional error-function-type integrals, based on
his work of random permutations. Br\'ezin and Hikami \cite{BH} apply
correlation functions of GUE ensemble to find explicit formulae of
$n$-point functions. Equation (44) in their paper agrees with the
$n=2$ case of our Theorem \ref{npoint2}.

Consider the following ``normalized'' $n$-point function
$$G(x_1,\dots,x_n)=\exp\left(\frac{-\sum_{j=1}^n
x_j^3}{24}\right) F(x_1,\dots,x_n).$$

The one-point function $G(x)=\frac{1}{x^2}$ is due to Witten, we
have also Dijkgraaf's two-point function
$$G(x,y)=\frac{1}{x+y}\sum_{k\geq0}\frac{k!}{(2k+1)!}\left(\frac{1}{2}xy(x+y)\right)^k$$
and Zagier's three-point function \cite{Za} which we learned from
Faber,
\begin{equation*}
G(x,y,z)=\sum_{r,s\geq
0}\frac{r!S_r(x,y,z)}{4^r(2r+1)!!\cdot2}\cdot\frac{\Delta^s}{8^s(r+s+1)!},
\end{equation*}
where $S_r(x,y,z)$ and $\Delta$ are the homogeneous symmetric polynomials defined
by
\begin{align*}
S_r(x,y,z)&=\frac{(xy)^r(x+y)^{r+1}+(yz)^r(y+z)^{r+1}+(zx)^r(z+x)^{r+1}}{x+y+z}\in\mathbb Z[x,y,z],\\
\Delta(x,y,z)&=(x+y)(y+z)(z+x)=\frac{(x+y+z)^3}{3}-\frac{x^3+y^3+z^3}{3}.
\end{align*}

The two and three point functions are discovered in the early 1990's.
Faber \cite{Fa} pioneered their use in the intersection theory of
moduli spaces of curves.

By studying Witten's KdV coefficient equation, regarded as an
ordinary differential equation, we get a recursive formula for
normalized $n$-point functions.

\begin{theorem} \label{npoint} For $n\geq2$,
\begin{equation*}
G(x_1,\dots,x_n)=\sum_{r,s\geq0}\frac{(2r+n-3)!!}{4^s(2r+2s+n-1)!!}P_r(x_1,\dots,x_n)\Delta(x_1,\dots,x_n)^s,
\end{equation*}
where $P_r$ and $\Delta$ are homogeneous symmetric polynomials
defined by
\begin{align*}
\Delta(x_1,\dots,x_n)&=\frac{(\sum_{j=1}^nx_j)^3-\sum_{j=1}^nx_j^3}{3},\nonumber\\
P_r(x_1,\dots,x_n)&=\left(\frac{1}{2\sum_{j=1}^nx_j}\sum_{\underline{n}=I\coprod
J}(\sum_{i\in I}x_i)^2(\sum_{i\in J}x_i)^2G(x_I)
G(x_J)\right)_{3r+n-3}\nonumber\\
&=\frac{1}{2\sum_{j=1}^nx_j}\sum_{\underline{n}=I\coprod
J}(\sum_{i\in I}x_i)^2(\sum_{i\in J}x_i)^2\sum_{r'=0}^r
G_{r'}(x_I)G_{r-r'}(x_J),
\end{align*}
where $I,J\ne\emptyset$, $\underline{n}=\{1,2,\ldots,n\}$ and $G_g(x_I)$ denotes the degree $3g+|I|-3$ homogeneous
component of the normalized $|I|$-point function
$G(x_{k_1},\dots,x_{k_{|I|}})$, where $k_j\in I$.
\end{theorem}

Note that the degree $3r+n-3$ polynomial $P_r(x_1,\dots,x_n)$ is
expressible by normalized $|I|$-point functions $G(x_I)$ with
$|I|<n$. So we can recursively obtain an explicit formula of the
$n$-point function
$$F(x_1,\dots,x_n)=\exp\left(\frac{\sum_{j=1}^n
x_j^3}{24}\right) G(x_1,\dots,x_n),$$ thus we have an elementary
algorithm to calculate all intersection numbers of $\psi$ classes.

Since $P_0(x,y)=\frac{1}{x+y}, P_r(x,y)=0$ for $r>0$, we get
Dijkgraaf's $2$-point function. From
$$P_r(x,y,z)=\frac{r!}{2^r(2r+1)!}\cdot\frac{(xy)^r(x+y)^{r+1}+(yz)^r(y+z)^{r+1}+(zx)^r(z+x)^{r+1}}{x+y+z},$$
we also get Zagier's $3$-point function.

We point out that the above recursive formula of normalized
$n$-point functions is essentially equivalent to the first
(classical) KdV equation \eqref{KdV} in Witten-Kontsevich theorem. See the
discussion at the latter part of Section 2 in \cite{LiuXu3}.

The results of this paper have applications to the tautological ring of moduli spaces of curves, Hodge integrals
and Gromov-Witten theory.

We will give a proof of Theorem \ref{npoint} in Section 2. Sections
3 contains some new identities of the intersection numbers of the
$\psi$ classes derived from the $n$-point functions. In Section 4,
we give a combinatorial interpretation of $n$-point functions in
terms of summation over binary trees. In section 5, we prove an
effective recursion formula for computing integrals of $\psi$
classes. In Section 6, we propose some conjectural generalization of
our results to Gromov-Witten invariants and Witten's $r$-spin
intersection
numbers. These conjectures have been proved recently by X. Liu and Pandharipande \\

\noindent{\bf Acknowledgements.} The authors would like to thank
Professors Sergei Lando, Jun Li, Chiu-Chu Melissa Liu, Xiaobo Liu,
Ravi Vakil and Jian Zhou for helpful communications. We also thank
Professors Edward Witten, Don Zagier for their interests in this
work and Professor Carel Faber for communicating Zagier's
three-point function to us.

\vskip 30pt
\section{Recursive formulae of n-point functions}

Theorem \ref{npoint} has several equivalent formulations.

\begin{proposition} \label{nptequiv} Let $n\geq2$. Then the recursion relation in Theorem \ref{npoint} is equivalent to either one of the following
statements.
\begin{enumerate}
\item[i)] The normalized $n$-point functions satisfy the following recursion relation
$$G_g(x_1,\dots,x_n)=
\frac{1}{(2g+n-1)}P_g(x_1,\dots,x_n)+\frac{\Delta(x_1,\dots,x_n)}{4(2g+n-1)}G_{g-1}(x_1,\dots,x_n).$$

\item[ii)] The $n$-point functions $F_g(x_1,\dots,x_n)$ satisfy the following recursion relation
\begin{multline*}
 (2g+n-1)\left(\sum_{i=1}^n x_i\right)F_g(x_1,\dots,x_n)=\frac{1}{12}\left(\sum_{i=1}^n x_i\right)^4F_{g-1}(x_1,\dots,x_n)\\
+\frac{1}{2}\sum_{g'=0}^g\sum_{\underline{n}=I\coprod J} \left(\sum_{i\in I}
x_i\right)^2\left(\sum_{i\in J} x_i\right)^2 F_{g'}(x_I)F_{g-g'}(x_J).
\end{multline*}
\end{enumerate}
\end{proposition}
\begin{proof}
For Theorem \ref{npoint} $\Rightarrow$ (i), we have
\begin{multline*}
G_g(x_1,\dots,x_n)=\sum_{r+s=g}\frac{(2r+n-3)!!}{4^s(2g+n-1)!!}P_r(x_1,\dots,x_n)\Delta(x_1,\dots,x_n)^s\\
=\frac{1}{2g+n-1}P_g(x_1,\dots,x_n)+\sum_{r+s=g-1}\frac{(2r+n-3)!!}{4^{s+1}(2g+n-1)!!}P_r(x_1,\dots,x_n)\Delta(x_1,\dots,x_n)^{s+1}\\
=\frac{1}{(2g+n-1)}P_g(x_1,\dots,x_n)+\frac{\Delta(x_1,\dots,x_n)}{4(2g+n-1)}G_{g-1}(x_1,\dots,x_n).
\end{multline*}

The proof that (i)
implies Theorem \ref{npoint} is also easy.

The equivalence of (i) and (ii) is the Proposition 2.3 of \cite{LiuXu3}.
\end{proof}

\begin{corollary} \label{nptequiv2}
For $n\geq2$,
\begin{equation*}
F(x_1,\dots,x_n)=\sum_{r,s\geq0}\frac{(2r+n-3)!!}{12^s(2r+2s+n-1)!!}S_r(x_1,\dots,x_n)\left(\sum_{j=1}^n
x_j\right)^{3s},
\end{equation*}
where $S_r$ is a homogeneous symmetric polynomial defined by
\begin{align*}
S_r(x_1,\dots,x_n)&=\left(\frac{1}{2\sum_{j=1}^nx_j}\sum_{\underline{n}=I\coprod
J}(\sum_{i\in I}x_i)^2 (\sum_{i\in J}x_i)^2 F(x_I)
F(x_J)\right)_{3r+n-3}\\
&=\frac{1}{2\sum_{j=1}^nx_j}\sum_{\underline{n}=I\coprod
J}(\sum_{i\in I}x_i)^2(\sum_{i\in J}x_i)^2\sum_{r'=0}^r
F_{r'}(x_I)F_{r-r'}(x_J),
\end{align*}
where $I,J\ne\emptyset$.
\end{corollary}
\begin{proof}
This follows directly from Proposition \ref{nptequiv} (ii).
\end{proof}

\begin{corollary}  We have
\begin{equation*}
\sum_{n\geq1}\sum_{\underline{n}=I\coprod J} \left(\sum_{i\in J} x_i\right)^4 F(-(x_1+\cdots+x_n),x_I)F(x_J)=1.
\end{equation*}
\end{corollary}
\begin{proof} Note that Proposition \ref{nptequiv} (ii) implies that for $2g+n-1>0$,
$$\sum_{g'=0}^g\sum_{\underline{n}=I\coprod J} \left(\sum_{i\in J} x_i\right)^4 F_{g'}(-(x_1+\cdots+x_n),x_I)F_{g-g'}(x_J)=0.$$
The right hand side $1$ comes from the case $n=1, g=0$.
\end{proof}

Recall that KdV hierarchy is captured in Witten's KdV coefficient equation
(see \cite{FaPa, Wi})
\begin{multline*}
(2d_1+1)\langle\tau_{d_1}\tau_0^2\prod_{j=2}^n\tau_{d_j}\rangle =
{\frac14}\langle\tau_{d_1-1}\tau_0^4\prod_{j=2}^n\tau_{d_j}\rangle\nonumber\\
+\sum_{\{2,\dots,n\}=I\coprod
J} \left(\langle\tau_{d_1-1}\tau_0\prod_{i\in
I}\tau_{d_i}\rangle\langle\tau_0^3\prod_{i\in J}\tau_{d_i}\rangle +
 2\langle\tau_{d_1-1}\tau_0^2\prod_{i\in I}\tau_{d_i}\rangle\langle\tau_0^2\prod_{i\in J}\tau_{d_i}\rangle\right),
\end{multline*}
which is equivalent to the following differential equation of
$n$-point functions (regarded as an ODE in $y$).
\begin{multline} \label{ODE}
y\frac{\partial}{\partial
y}\left((y+\sum_{j=1}^{n}x_j)^2F_g(y,x_1,\dots,x_n)\right)\\
=\frac{y}{8}(y+\sum_{j=1}^{n}x_j)^{4}F_{g-1}(y,x_1,\dots,x_n)+\frac{y}{2}(y+\sum_{j=1}^{n}x_j)F_g(y,x_1,\dots,x_n)\\
+\frac{y}{2}\sum_{\underline{n}=I\coprod J}\left(\left(y+\sum_{i\in
I}x_i\right)\left(\sum_{i\in J}x_i\right)^3+2\left(y+\sum_{i\in
I}x_i\right)^2\left(\sum_{i\in
J}x_i\right)^2\right)F_{g'}(y,x_I)F_{g-g'}(x_J)\\-\frac{1}{2}\left(y+\sum_{j=1}^{n}x_j\right)^2F_g(y,x_1,\dots,x_n)
\end{multline}

\subsection{Proof of Theorem \ref{npoint}} By Proposition
\ref{nptequiv}, in order to prove Theorem \ref{npoint}, it is
sufficient to verify that $F(x_1,\dots,x_n)$, as recursively defined
in Proposition \ref{nptequiv} (ii), satisfies the above differential
equation. The verification is tedious but straightforward. The
details are in the appendix.

Moreover, we need to check the initial value condition (the string equation)
$$F(x_1,\dots,x_n,0)=(\sum_{j=1}^n x_j) F(x_1,\dots,x_n).$$

By induction, we have
\begin{multline*}
\left(\sum_{j=1}^n
x_j\right)F_g(x_1,\dots,x_n,0)=\frac{1}{2g+n}\left(\frac{\left(\sum_{j=1}^n
x_j\right)^4}{12}
F_{g-1}(x_1,\dots,x_n,0)\right.\\+\left(\sum_{j=1}^n x_j\right)^2
F_g(x_1,\dots,x_n)
\left.+\frac{1}{2}\sum_{h=0}^g\sum_{\underline{n}=I\coprod J}
\left(\sum_{i\in I} x_i\right)^2\left(\sum_{i\in J} x_i\right)^2
F_h(x_I,0)F_{g-h}(x_J)\right.\\\left.+\frac{1}{2}\sum_{h=0}^g\sum_{\underline{n}=I\coprod
J} \left(\sum_{i\in I} x_i\right)^2\left(\sum_{i\in J} x_i\right)^2
F_h(x_I)F_{g-h}(x_J,0)\right)\\
=\frac{1}{2g+n}\left(\left(\sum_{j=1}^n x_j\right)^2
F_g(x_1,\dots,x_n)+(2g+n-1)\left(\sum_{j=1}^n x_j\right)^2
F_g(x_1,\dots,x_n)\right) \\=\left(\sum_{j=1}^n x_j\right)^2
F_g(x_1,\dots,x_n).
\end{multline*}

By the uniqueness of ODE solutions, we have proved Theorem \ref{npoint}.
\smallskip

In the meantime, we also proved the following result, which explains why in order to prove
the Witten-Kontsevich theorem, it suffices to prove that the generating function \eqref{gen} satisfies the classical KdV equation \eqref{KdV}, as was done in \cite{KaL}.

\begin{corollary}
Under constraints of the string and dilaton equations,
$$\left(-\frac{\partial}{\partial t_0}+\sum_{i=0}^{\infty} t_{i+1}+\frac{t_0^2}{2}\right) \exp F(t_i)=0$$
$$\left(-\frac{3}{2}\frac{\partial}{\partial t_1}+\sum_{i=0}^{\infty}\frac{2i+1}{2}t_i\frac{\partial}{\partial t_i}+\frac{1}{16}\right)\exp F(t_i)=0,$$
any quasi-homogeneous solution $F(t_i)=\sum_{g=0}^{\infty}F_g(t_i)$ to the classical KdV equation
automatically satisfies the whole KdV hierarchy.
\end{corollary}

There is another slightly different formula of $n$-point functions.
When $n=3$, this has also been obtained by Zagier \cite{Za}.
\begin{theorem} \label{npoint2} For $n\geq2$,
\begin{equation*}
F(x_1,\dots,x_n)=\exp\left(\frac{(\sum_{j=1}^n
x_j)^3}{24}\right)\sum_{r,s\geq0}\frac{(-1)^s
P_r(x_1,\dots,x_n)\Delta(x_1,\dots,x_n)^s}{8^s(2r+2s+n-1)s!},
\end{equation*}
where $P_r$ and $\Delta$ are the same polynomials as defined in
Theorem \ref{npoint}.
\end{theorem}

It is easy to see that Theorem  \ref{npoint2} follows from Theorem  \ref{npoint} and the following lemma.
\begin{lemma} Let $n\geq2$ and $r,s\geq0$. Then the following identity holds,
\begin{equation*}
\frac{(-1)^s}{8^s(2r+2s+n-1)s!}=\sum_{k=0}^s\frac{(-1)^k}{8^k
k!}\cdot \frac{(2r+n-3)!!}{4^{s-k}(2r+2s-2k+n-1)!!}
\end{equation*}
\end{lemma}
\begin{proof}
Let $p=2r+n\geq2$ and
$$f(p,s)=\sum_{k=0}^s\frac{(-1)^k}{2^k k!(p+2s-2k-1)!!}.$$
We have
\begin{align*}
f(p,s)&=\sum_{k=0}^s\frac{(-1)^k(p+2s+1)}{2^k k!(p+2s-2k+1)!!}+\sum_{k=0}^s\frac{2k(-1)^{k-1} }{2^k k!(p+2s-2k+1)!!}\\
&=(p+2s+1)\left(f(p,s+1)-\frac{(-1)^{s+1}}{2^{s+1}(s+1)!(p-1)!!}\right)+f(p,s)-\frac{(-1)^s}{2^s
s!(p-1)!!}.
\end{align*}
So we have the following identity
$$f(p,s+1)=\frac{(-1)^{s+1}}{2^{s+1}(p+2s+1)(s+1)!(p-3)!!},$$
which is just the identity we want if $s+1$ is replaced by $s$.
\end{proof}

\vskip 30pt
\section{New properties of the n-point functions}

In this section we derive various new identities about the
coefficients of the $n$-point functions. An important application is
a proof of the famous Faber intersection number conjecture
\cite{Fa}. Recently, Zhou \cite{Zh} used our results on $n$-point
functions in his computation of Hurwitz-Hodge integrals.

Let $\mathcal{C} \left(\prod_{j=1}^n
x_j^{d_j},p(x_1,\dots,x_n)\right)$ denote the coefficient of
$\prod_{j=1}^n x_j^{d_j}$ in a polynomial or formal power series
$p(x_1,\dots,x_n)$. From the inductive structure in the definition
of $n$-point functions, we have the following basic properties of
$n$-point functions.

First consider the normalized $(n+1)$-point function
$G(y,x_1,\dots,x_n)$. Here we use $y$ to denote a distinguished point.

\begin{theorem} \label{coeff} Let $2g-2+n\geq0$.
\begin{enumerate}
\item[i)] Let $k>2g-2+n$, $d_j\geq0$ and $\sum_{j=1}^n
d_j=3g-2+n-k$. Then
\begin{align*}
\mathcal C\left(y^k\prod_{j=1}^{n}x_j^{d_j}, G_g(y,x_1,\dots,x_n)\right)&=0,\\
\mathcal C\left(y^k\prod_{j=1}^{n}x_j^{d_j}, P_g(y,x_1,\dots,x_n)\right)&=0.
\end{align*}

\item[ii)] Let $d_j\geq0$, $\sum_{j=1}^n
d_j=g$ and $a=\#\{j\mid d_j=0\}$. Then
\begin{align*}
\mathcal C\left(y^{2g-2+n}\prod_{j=1}^{n}x_j^{d_j},
G_g(y,x_1,\dots,x_n)\right)&=\frac{1}{4^g\prod_{j=1}^n(2d_j+1)!!},\\
\mathcal C\left(y^{2g-2+n}\prod_{j=1}^{n}x_j^{d_j},
P_g(y,x_1,\dots,x_n)\right)&=\frac{a}{4^g\prod_{j=1}^n(2d_j+1)!!}.
\end{align*}

\item[iii)] Let $d_j\geq0$, $\sum_{j=1}^n
d_j=g+1$, $a=\#\{j\mid d_j=0\}$ and $b=\#\{j\mid d_j=1\}$. Then
\begin{align*}
\mathcal C\left(y^{2g-3+n}\prod_{j=1}^{n}x_j^{d_j},
G_g(y,x_1,\dots,x_n)\right)&=\frac{2g^2+(2n-1)g+\frac{n^2-n}{2}-3+\frac{5a-a^2}{2}}{4^g\prod_{j=1}^n(2d_j+1)!!},\\
\mathcal C\left(y^{2g-3+n}\prod_{j=1}^{n}x_j^{d_j},
P_g(y,x_1,\dots,x_n)\right)&=\frac{a(2g^2+2ng-g+\frac{n^2-n-a^2+5a}2+3b-3)-3b}{4^g\prod_{j=1}^n(2d_j+1)!!}.
\end{align*}
\end{enumerate}
\end{theorem}
\begin{proof}
The proof uses Proposition \ref{nptequiv} (i) and proceeds by induction on $g$ and $n$. Note that
$$\Delta(y, x_1,\dots,x_n)=y^2(\sum_{j=1}^n x_j)+y(\sum_{j=1}^n x_j)^2+\Delta(x_1,\dots,x_n).$$
The vanishing identities (i) are obvious. We now prove (ii) inductively.
\begin{align*}
\mathcal C\left(y^{2g-2+n}\prod_{j=1}^{n}x_j^{d_j},
P_g(y,x_1,\dots,x_n)\right)&=\sum_{j=1}^n\mathcal C\left(y^{2g-2+n}\prod_{j=1}^{n}x_j^{d_j}, G_g(y,x_1,\dots,\hat{x_j},\dots,x_n)\right)\\
&=\frac{a}{4^g\prod_{j=1}^n(2d_j+1)!!},
\end{align*}
where $a=\#\{j\mid d_j=0\}$.

\begin{align*}
\mathcal C&\left(y^{2g-2+n}, G_g(y,x_1,\dots,x_n)\right)\\
&=\sum_{r+s=g}\frac{(2r+n-2)!!}{4^s(2g+n)!!}\sum_{\sum d_j=r}\frac{a\cdot\prod_{j=1}^n x_j^{d_j}}{4^r\prod_{j=1}^n(2d_j+1)!!}(\sum_{j=1}^nx_j)^s\\
&=\frac{1}{2g+n}\sum_{\sum d_j=g}\frac{a\cdot\prod_{j=1}^n x_j^{d_j}}{4^g\prod_{j=1}^n(2d_j+1)!!}+\frac{\sum_{j=1}^n x_j}{4(2g+n)}\sum_{\sum d_j=g-1}\frac{\prod_{j=1}^n x_j^{d_j}}{4^{g-1}\prod_{j=1}^n(2d_j+1)!!}\\
&=\frac{1}{2g+n}\left(\sum_{\sum d_j=g}\frac{a\cdot\prod_{j=1}^n x_j^{d_j}}{4^g\prod_{j=1}^n(2d_j+1)!!}+\sum_{\sum d_j=g}\frac{(2g+n-a)\prod_{j=1}^n x_j^{d_j}}{4^g\prod_{j=1}^n(2d_j+1)!!}\right)\\
&=\sum_{\sum d_j=g}\frac{\prod_{j=1}^n x_j^{d_j}}{4^g\prod_{j=1}^n(2d_j+1)!!}.
\end{align*}

The identities (iii) can be proved similarly.
\end{proof}

\begin{corollary}
Let $2g-2+n\geq0$.
\begin{enumerate}
\item[i)] Let $k>2g-2+n$, $d_j\geq0$ and $\sum_{j=1}^n
d_j=3g-2+n-k$. Then
\begin{equation*}
\sum_{r=0}^g \frac{(-1)^r}{24^r r!}\langle\tau_0^{3r}\tau_k\tau_{d_1}\cdots\tau_{d_n}\rangle_{g-r}=0.
\end{equation*}

\item[ii)] Let $d_j\geq0$ and $\sum_{j=1}^n
d_j=g$. Then
\begin{equation*}
\sum_{r=0}^g \frac{(-1)^r}{24^r r!}\langle\tau_0^{3r}\tau_{2g-2+n}\tau_{d_1}\cdots\tau_{d_n}\rangle_{g-r}=\frac{(-1)^g}{8^g\prod_{j=1}^n d_j ! \cdot(2d_j+1)}.
\end{equation*}
\end{enumerate}
\end{corollary}

\begin{proof} We have 
$$\exp\left(\frac{-(y+\sum_{j=1}^n
x_j)^3}{24}\right)F(y, x_1,\dots,x_n)=\exp\left(\frac{-\Delta(y, x_1,\dots,x_n)}{8}\right)G(y, x_1,\dots,x_n).$$
We need to extract coefficients from both sides and the corollary follows by an induction using Theorem \ref{coeff}.
\end{proof}

We may regard $F(y,x_1,\dots,x_n)$ and $G(y,x_1,\dots,x_n)$ as formal
series in $\mathbb Q[x_1,\dots,x_n][[y,y^{-1}]]$ with $\deg y<\infty$.
In particular,
\begin{equation} \label{con}
F_0(y)=G_0(y)=\frac{1}{y^2},\qquad F_{0}(x,y)=G_0(x,y)=\frac{1}{x+y}=\sum_{k=0}^\infty (-1)^k\frac{x^k}{y^{k+1}}.
\end{equation}

We can again use Proposition \ref{nptequiv} to prove the following proposition inductively, which is crucial in our proof of the famous Faber intersection number conjecture \cite{LiuXu3}.
\begin{proposition} \label{nptcoeff}  Let $a,b\in\mathbb Z$.
\begin{enumerate}

\item[i)]
Let  $k\geq 2g-3+a+b$. Then
$$
\mathcal C\left(y^k,\quad \sum_{g'=0}^g\sum_{\underline{n}=I\coprod
J} (y+\sum_{i\in I}x_i)^a (-y+\sum_{i\in J}x_i)^b
F_{g'}(y,x_I)F_{g-g'}(-y,x_J)\right)=0.$$
$$
\mathcal C\left(y^k,\quad \sum_{g'=0}^g\sum_{\underline{n}=I\coprod
J}(y+\sum_{i\in I}x_i)^a(-y+\sum_{i\in J}x_i)^b G_{g'}(y,x_I)
G_{g-g'}(-y,x_J)\right)=0.$$

\item[ii)]
Let $d_j\geq 1$ and $\sum_{j}d_j=g+n$. Then
\begin{multline*}
\mathcal C\left(y^{2g-4+a+b}\prod_{j=1}^nx_j^{d_j},
\sum_{g'=0}^g\sum_{\underline{n}=I\coprod J} (y+\sum_{i\in I}x_i)^a
(-y+\sum_{i\in J}x_i)^b
F_{g'}(y,x_I)F_{g-g'}(-y,x_J)\right)\\
=\mathcal C\left(y^{2g-4+a+b}\prod_{j=1}^nx_j^{d_j},
\sum_{g'=0}^g\sum_{\underline{n}=I\coprod J} (y+\sum_{i\in I}x_i)^a
(-y+\sum_{i\in J}x_i)^b G_{g'}(y,x_I)G_{g-g'}(-y,x_J)\right)\\=\frac{(-1)^b(2g-3+n+a+b)!}{4^g(2g-3+a+b)!\prod_{j=1}^n(2d_j-1)!!}.
\end{multline*}
\end{enumerate}
\end{proposition}

When $n=2$, Proposition \ref{nptcoeff} can be checked directly. For example, take
$a=b=2$,
\begin{align*}
&\mathcal C\left(y^{2g+2}, \sum_{\underline{2}=I\coprod
J}(y+\sum_{i\in I}x_i)^2 (-y+\sum_{i\in J}x_i)^2 \sum_{g'=0}^g
G_{g'}(y,x_I)G_{g-g'}(-y,x_J)\right)\\
&=2\sum_{r+s=g}\left(\frac{(2r)!!}{4^s(2g+2)!!}\frac{1}{4^r(2r+1)!!}(x_1^r+x_2^r)(x_1+x_2)^s-\frac{1}{4^r(2r+1)!!}\frac{1}{4^s(2s+1)!!}x_1^r
x_2^s\right)\\
&=0.
\end{align*}

We can extract coefficients of $n$-point functions to get identities for intersection numbers of $\psi$ classes.
A detailed discussion can be found in \cite{LiuXu3}. We record two such identities here.
\begin{corollary} We have
\begin{enumerate}

\item[i)] Let $d_j\geq0$, $\#\{j\mid d_j=0\}\leq 1$ and $\sum_{j=1}^{n}(d_j-1)=g-1$. Then
$$\sum_{j=0}^{2g}(-1)^j\langle\tau_{2g-j}\tau_j\prod_{i=1}^n\tau_{d_i}\rangle_{g}=
\frac{(2g+n-1)!}{4^{g}(2g+1)!\prod_{j=1}^n(2d_j-1)!!}.$$ If
$\#\{j\mid d_j=0\}=2$ and $a=\#\{j\mid d_j=1\}$, then the right hand
side becomes
$$\frac{(2g+n-1)!}{4^g(2g+1)!\prod_{j=1}^n(2d_j-1)!!}\cdot\frac{2g+n-a}{2g+n-1-a}.$$

\item[ii)]
Let $d_j\geq1$ and $\sum_{j=1}^{n}(d_j-1)=g$. Then
\begin{multline*}
\frac{(2g-3+n)!}{2^{2g+1}(2g-3)!\prod_{j=1}^n(2d_j-1)!!}
=\langle\tau_{2g-2}\prod_{j=1}^n\tau_{d_j}\rangle_g-\sum_{j=1}^{n}
\langle\tau_{d_j+2g-3}\prod_{i\neq j}\tau_{d_i}\rangle_g\\
+\frac{1}{2}\sum_{\underline{n}=I\coprod
J}\sum_{j=0}^{2g-4}(-1)^j\langle\tau_{j}\prod_{i\in
I}\tau_{d_i}\rangle_{g'}\langle\tau_{2g-4-j}\prod_{i\in
J}\tau_{d_i}\rangle_{g-g'}.
\end{multline*}

If $\#\{j\mid d_j=0\}=1$ and $a=\#\{j\mid
d_j=1\}$, then the left hand side becomes
$$\frac{(2g-3+n)!}{2^{2g+1}(2g-3)!\prod_{j=1}^n(2d_j-1)!!}\cdot\frac{2g+n+1-a}{2g+n-3-a}.
$$
\end{enumerate}
\end{corollary}

\vskip 30pt
\section{$n$-point function as summation over binary trees}

Recall that in graph theory, a ``tree'' is defined to be a graph
without cycles. A ``binary tree'' $T$ is a tree such that
each node $v\in V(T)$ either has no children ($v\in L(T)$ is a leaf)
or has two children ($v\notin L(T)$), so we must have
$|V(T)|=2|L(T)|-1$.

Denote by $r_T$ the unique root of $T$. For each $v\in V(T)$, define $D(v)\subset V(T)$ to be the set of all
descendants of $v$ and define $L(v)=D(v)\cap L(T)$. In particular,
if $v$ is a leaf, then $D(v)=L(v)=\{v\}$; if $v=r_T$, then $D(v)=V(T)$ and $L(v)=L(T)$.

\begin{definition}
Let $T$ be a binary tree. Let $n=|L(T)|$ be the number of leaves. We
assign an integer $g(v)\geq0$ to each node $v\in V(T)$ and label the
$n$ leaves with distinct values $\ell(v)\in\{1,\dots,n\}$. Then we
call such $T$ a ``weighted marked binary tree'' (abbreviated ``WMB
tree'') and call $g(T)=\sum_{v\in V(T)}g(v)$ the total weight of
$T$.
\end{definition}

Now we can state our main result in this section.
\begin{theorem} \label{tree}
Let $g\geq0, n\geq1$. Denote by WMB$(g,n)$ the set of isomorphism
classes of all WMB trees with total weight $g$ and $n$ leaves. Then
we have the following expression of $n$-point functions:
\begin{multline*}
12^g\left(\prod_{j=1}^n x_j\right)\cdot(x_1+\dots+x_n)^2
F_g(x_1,\dots,x_n)\\= \sum_{T\in\ {\rm WMB}(g,n)} \prod_{v\in
V(T)}\frac{\left(|L(v)|-3+\sum\limits_{\substack{w\in D(v)\\w\neq
v}}2g(w)\right)!!}{\left(|L(v)|-1+\sum\limits_{w\in
D(v)}2g(w)\right)!!}\left(\sum_{w\in
L(v)}x_{\ell(w)}\right)^{3g(v)+1}.
\end{multline*}
Note that $(-2)!!=(-1)!!=0!!=1$ by definition.
\end{theorem}
\begin{proof}
When $n=1$, the identity holds obviously.

By noting that bipartition of indices corresponds to siblings in
binary trees and applying corollary \ref{nptequiv2} recursively, we
may get
\begin{multline} \label{comb}
(x_1+\dots+x_n)^2 F_g(x_1,\dots,x_n)\\= \sum_{T\in\ {\rm WMB}(g,n)}
\prod_{v\in
L(T)}\frac{x_{\ell(v)}^{3g(v)}}{24^{g(v)}g(v)!}\prod_{v\notin
L(T)}\frac{\left(|L(v)|-3+\sum\limits_{\substack{w\in D(v)\\w\neq
v}}2g(w)\right)!!\left(\sum\limits_{w\in
L(v)}x_{\ell(w)}\right)^{3g(v)+1}}{12^{g(v)}\left(|L(v)|-1+\sum\limits_{w\in
D(v)}2g(w)\right)!!}\\
=\frac{1}{12^g}\sum_{T\in\ {\rm WMB}(g,n)} \prod_{v\in
L(T)}\frac{x_{\ell(v)}^{3g(v)}}{(2g(v))!!}\prod_{v\notin
L(T)}\frac{\left(|L(v)|-3+\sum\limits_{\substack{w\in D(v)\\w\neq
v}}2g(w)\right)!!\left(\sum\limits_{w\in
L(v)}x_{\ell(w)}\right)^{3g(v)+1}}{\left(|L(v)|-1+\sum\limits_{w\in
D(v)}2g(w)\right)!!}.
\end{multline}

So we get the desired identity. The details are left to the
interested readers.
\end{proof}

We now illustrate the above theorem by two examples. We will compute
the right hand side of the slightly simpler identity \eqref{comb},
which avoids the factor $\prod_{j=1}^n x_j$.
\begin{example}

Take $(g,n)=(1,3)$. In the following WMB trees, the number in the
circle denotes the label of a leaf, while the number beside each
node represents its weight.

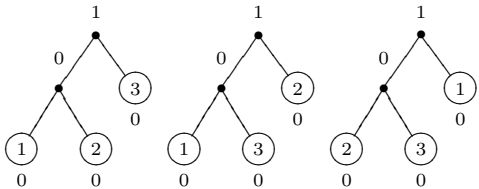
\begin{minipage}{0.8in}
\begin{displaymath}
  \vcenter{\tiny\xymatrix@C=1mm@R=0mm{
  & &1& &\\
  &  &  *={\bullet}\ar@{-}[ddl]\ar@{-}[ddr] & & \\
  &0& & &\\
  &  *={\bullet}\ar@{-}[ddl]\ar@{-}[ddr] & &
  *=<12pt>[o][F-]{3}&  \\
  & & &0&\\
  *=<12pt>[o][F-]{1} & &  *=<12pt>[o][F-]{2}
& &\\
0& &0& &
         }}\,
\end{displaymath}
\end{minipage}
\begin{minipage}{0.8in}
\begin{displaymath}
  \vcenter{\tiny\xymatrix@C=1mm@R=0mm{
  & &1& &\\
  &  &  *={\bullet}\ar@{-}[ddl]\ar@{-}[ddr] & & \\
  &0& & &\\
  &  *={\bullet}\ar@{-}[ddl]\ar@{-}[ddr] & &
  *=<12pt>[o][F-]{2}&  \\
  & & &0&\\
  *=<12pt>[o][F-]{1} & &  *=<12pt>[o][F-]{3}
& &\\
0& &0& &
         }}\,
\end{displaymath}
\end{minipage}
\begin{minipage}{0.8in}

\begin{displaymath}
  \vcenter{\tiny\xymatrix@C=1mm@R=0mm{
  & &1& &\\
  &  &  *={\bullet}\ar@{-}[ddl]\ar@{-}[ddr] & & \\
  &0& & &\\
  &  *={\bullet}\ar@{-}[ddl]\ar@{-}[ddr] & &
  *=<12pt>[o][F-]{1}&  \\
  & & &0&\\
  *=<12pt>[o][F-]{2} & &  *=<12pt>[o][F-]{3}
& &\\
0& &0& &
         }}\,
\end{displaymath}
\end{minipage}
\smallskip

$I_1=(x_1+x_2+x_3)^4((x_1+x_2)+(x_1+x_3)+(x_2+x_3))/96.$

\begin{minipage}{0.8in}

\begin{displaymath}
  \vcenter{\tiny\xymatrix@C=1mm@R=0mm{
  & &0& &\\
  &  &  *={\bullet}\ar@{-}[ddl]\ar@{-}[ddr] & & \\
  &1& & &\\
  &  *={\bullet}\ar@{-}[ddl]\ar@{-}[ddr] & &
  *=<12pt>[o][F-]{3}&  \\
  & & &0&\\
  *=<12pt>[o][F-]{1} & &  *=<12pt>[o][F-]{2}
& &\\
0& &0& &
         }}\,
\end{displaymath}
\end{minipage}
\begin{minipage}{0.8in}

\begin{displaymath}
  \vcenter{\tiny\xymatrix@C=1mm@R=0mm{
  & &0& &\\
  &  &  *={\bullet}\ar@{-}[ddl]\ar@{-}[ddr] & & \\
  &1& & &\\
  &  *={\bullet}\ar@{-}[ddl]\ar@{-}[ddr] & &
  *=<12pt>[o][F-]{2}&  \\
  & & &0&\\
  *=<12pt>[o][F-]{1} & &  *=<12pt>[o][F-]{3}
& &\\
0& &0& &
         }}\,
\end{displaymath}
\end{minipage}
\begin{minipage}{0.8in}

\begin{displaymath}
  \vcenter{\tiny\xymatrix@C=1mm@R=0mm{
  & &0& &\\
  &  &  *={\bullet}\ar@{-}[ddl]\ar@{-}[ddr] & & \\
  &1& & &\\
  &  *={\bullet}\ar@{-}[ddl]\ar@{-}[ddr] & &
  *=<12pt>[o][F-]{1}&  \\
  & & &0&\\
  *=<12pt>[o][F-]{2} & &  *=<12pt>[o][F-]{3}
& &\\
0& &0& &
         }}\,
\end{displaymath}
\end{minipage}
\smallskip

$I_2=(x_1+x_2+x_3)((x_1+x_2)^4+(x_1+x_3)^4+(x_2+x_3)^4)/144.$

\begin{minipage}{0.8in}

\begin{displaymath}
  \vcenter{\tiny\xymatrix@C=1mm@R=0mm{
  & &0& &\\
  &  &  *={\bullet}\ar@{-}[ddl]\ar@{-}[ddr] & & \\
  &0& & &\\
  &  *={\bullet}\ar@{-}[ddl]\ar@{-}[ddr] & &
  *=<12pt>[o][F-]{1}&  \\
  & & &1&\\
  *=<12pt>[o][F-]{2} & &  *=<12pt>[o][F-]{3}
& &\\
0& &0& &
         }}\,
\end{displaymath}
\end{minipage}
\begin{minipage}{0.8in}

\begin{displaymath}
  \vcenter{\tiny\xymatrix@C=1mm@R=0mm{
  & &0& &\\
  &  &  *={\bullet}\ar@{-}[ddl]\ar@{-}[ddr] & & \\
  &0& & &\\
  &  *={\bullet}\ar@{-}[ddl]\ar@{-}[ddr] & &
  *=<12pt>[o][F-]{2}&  \\
  & & &1&\\
  *=<12pt>[o][F-]{1} & &  *=<12pt>[o][F-]{3}
& &\\
0& &0& &
         }}\,
\end{displaymath}
\end{minipage}
\begin{minipage}{0.8in}

\begin{displaymath}
  \vcenter{\tiny\xymatrix@C=1mm@R=0mm{
  & &0& &\\
  &  &  *={\bullet}\ar@{-}[ddl]\ar@{-}[ddr] & & \\
  &0& & &\\
  &  *={\bullet}\ar@{-}[ddl]\ar@{-}[ddr] & &
  *=<12pt>[o][F-]{3}&  \\
  & & &1&\\
  *=<12pt>[o][F-]{1} & &  *=<12pt>[o][F-]{2}
& &\\
0& &0& &
         }}\,
\end{displaymath}
\end{minipage}
\smallskip

$I_3=(x_1+x_2+x_3)(x_1^3(x_2+x_3)+x_2^3(x_1+x_3)+x_3^3(x_1+x_2)^4)/96.$

\begin{minipage}{0.8in}

\begin{displaymath}
  \vcenter{\tiny\xymatrix@C=1mm@R=0mm{
  & &0& &\\
  &  &  *={\bullet}\ar@{-}[ddl]\ar@{-}[ddr] & & \\
  &0& & &\\
  &  *={\bullet}\ar@{-}[ddl]\ar@{-}[ddr] & &
  *=<12pt>[o][F-]{3}&  \\
  & & &0&\\
  *=<12pt>[o][F-]{1} & &  *=<12pt>[o][F-]{2}
& &\\
1& &0& &
         }}\,
\end{displaymath}
\end{minipage}
\begin{minipage}{0.8in}

\begin{displaymath}
  \vcenter{\tiny\xymatrix@C=1mm@R=0mm{
  & &0& &\\
  &  &  *={\bullet}\ar@{-}[ddl]\ar@{-}[ddr] & & \\
  &0& & &\\
  &  *={\bullet}\ar@{-}[ddl]\ar@{-}[ddr] & &
  *=<12pt>[o][F-]{2}&  \\
  & & &0&\\
  *=<12pt>[o][F-]{1} & &  *=<12pt>[o][F-]{3}
& &\\
1& &0& &
         }}\,
\end{displaymath}
\end{minipage}
\begin{minipage}{0.8in}

\begin{displaymath}
  \vcenter{\tiny\xymatrix@C=1mm@R=0mm{
  & &0& &\\
  &  &  *={\bullet}\ar@{-}[ddl]\ar@{-}[ddr] & & \\
  &0& & &\\
  &  *={\bullet}\ar@{-}[ddl]\ar@{-}[ddr] & &
  *=<12pt>[o][F-]{3}&  \\
  & & &0&\\
  *=<12pt>[o][F-]{2} & &  *=<12pt>[o][F-]{1}
& &\\
1& &0& &
         }}\,
\end{displaymath}
\end{minipage}
\begin{minipage}{0.8in}

\begin{displaymath}
  \vcenter{\tiny\xymatrix@C=1mm@R=0mm{
  & &1& &\\
  &  &  *={\bullet}\ar@{-}[ddl]\ar@{-}[ddr] & & \\
  &0& & &\\
  &  *={\bullet}\ar@{-}[ddl]\ar@{-}[ddr] & &
  *=<12pt>[o][F-]{1}&  \\
  & & &0&\\
  *=<12pt>[o][F-]{2} & &  *=<12pt>[o][F-]{3}
& &\\
1& &0& &
         }}\,
\end{displaymath}
\end{minipage}
\begin{minipage}{0.8in}

\begin{displaymath}
  \vcenter{\tiny\xymatrix@C=1mm@R=0mm{
  & &0& &\\
  &  &  *={\bullet}\ar@{-}[ddl]\ar@{-}[ddr] & & \\
  &0& & &\\
  &  *={\bullet}\ar@{-}[ddl]\ar@{-}[ddr] & &
  *=<12pt>[o][F-]{2}&  \\
  & & &0&\\
  *=<12pt>[o][F-]{3} & &  *=<12pt>[o][F-]{1}
& &\\
1& &0& &
         }}\,
\end{displaymath}
\end{minipage}
\begin{minipage}{0.8in}

\begin{displaymath}
  \vcenter{\tiny\xymatrix@C=1mm@R=0mm{
  & &1& &\\
  &  &  *={\bullet}\ar@{-}[ddl]\ar@{-}[ddr] & & \\
  &0& & &\\
  &  *={\bullet}\ar@{-}[ddl]\ar@{-}[ddr] & &
  *=<12pt>[o][F-]{1}&  \\
  & & &0&\\
  *=<12pt>[o][F-]{3} & &  *=<12pt>[o][F-]{2}
& &\\
1& &0& &
         }}\,
\end{displaymath}
\end{minipage}
\smallskip

$I_4=(x_1+x_2+x_3)(x_1^3(2x_1+x_2+x_3)+x_2^3(2x_2+x_1+x_3)+x_3^3(2x_3+x_1+x_2))/288.$

\bigskip

From the above computation, we get
\begin{multline*}
F_1(x_1,x_2,x_3)=\frac{I_1+I_2+I_3+I_4}{(x_1+x_2+x_3)^2}\\
=\frac{x_1^3+x_2^3+x_3^3}{24}+\frac{x_1^2 x_2+x_1^2 x_3+x_1
x_2^2+x_1 x_3^2+x_2^2 x_3+x_2 x_3^2}{12} +\frac{x_1 x_2 x_3}{12},
\end{multline*}
which is easily seen to be correct.
\end{example}

\begin{example} Take $(g,n)=(2,2)$.

\begin{minipage}{0.8in}
\begin{displaymath}
  \vcenter{\tiny\xymatrix@C=1mm@R=0mm{
   &2& \\
    &  *={\bullet}\ar@{-}[ddl]\ar@{-}[ddr] &  \\
  & & \\
 *=<12pt>[o][F-]{1}  &   &
  *=<12pt>[o][F-]{2}  \\
   0&&0
         }}\,
\end{displaymath}
\end{minipage}
\begin{minipage}{0.8in}
\begin{displaymath}
  \vcenter{\tiny\xymatrix@C=1mm@R=0mm{
   &1& \\
    &  *={\bullet}\ar@{-}[ddl]\ar@{-}[ddr] &  \\
  & & \\
 *=<12pt>[o][F-]{1}  &   &
  *=<12pt>[o][F-]{2}  \\
   1&&0
         }}\,
\end{displaymath}
\end{minipage}
\begin{minipage}{0.8in}

\begin{displaymath}
  \vcenter{\tiny\xymatrix@C=1mm@R=0mm{
   &1& \\
    &  *={\bullet}\ar@{-}[ddl]\ar@{-}[ddr] &  \\
  & & \\
 *=<12pt>[o][F-]{1}  &   &
  *=<12pt>[o][F-]{2}  \\
   0&&1
         }}\,
\end{displaymath}
\end{minipage}
\smallskip

$$
I_1=\frac{(x_1+x_2)^7}{12^2\cdot
15}+\frac{x_1^3(x_1+x_2)^4}{12^2\cdot 2\cdot 15}
+\frac{x_2^3(x_1+x_2)^4}{12^2\cdot 2 \cdot 15}.
$$

\begin{minipage}{0.8in}
\begin{displaymath}
  \vcenter{\tiny\xymatrix@C=1mm@R=0mm{
   &0& \\
    &  *={\bullet}\ar@{-}[ddl]\ar@{-}[ddr] &  \\
  & & \\
 *=<12pt>[o][F-]{1}  &   &
  *=<12pt>[o][F-]{2}  \\
   2&&0
         }}\,
\end{displaymath}
\end{minipage}
\begin{minipage}{0.8in}
\begin{displaymath}
  \vcenter{\tiny\xymatrix@C=1mm@R=0mm{
   &0& \\
    &  *={\bullet}\ar@{-}[ddl]\ar@{-}[ddr] &  \\
  & & \\
 *=<12pt>[o][F-]{1}  &   &
  *=<12pt>[o][F-]{2}  \\
   1&&1
         }}\,
\end{displaymath}
\end{minipage}
\begin{minipage}{0.8in}

\begin{displaymath}
  \vcenter{\tiny\xymatrix@C=1mm@R=0mm{
   &0& \\
    &  *={\bullet}\ar@{-}[ddl]\ar@{-}[ddr] &  \\
  & & \\
 *=<12pt>[o][F-]{1}  &   &
  *=<12pt>[o][F-]{2}  \\
   0&&2
         }}\,
\end{displaymath}
\end{minipage}
\smallskip

$$I_2=\frac{3!!}{12^2\cdot 8\cdot 5!!}x_1^6(x_1+x_2)+\frac{3!!}{12^2\cdot 4\cdot 5!!}x_1^3x_2^3(x_1+x_2)+
\frac{3!!}{12^2\cdot 8\cdot 5!!}x_2^6(x_1+x_2).$$

\bigskip

From the above computation, we get
\begin{equation*}
F_1(x_1,x_2,x_3)=\frac{I_1+I_2+I_3+I_4}{(x_1+x_2+x_3)^2}
=\frac{x_1^5+x_2^5}{1152}+\frac{x_1^4 x_2+x_1
x_2^4}{384}+\frac{29x_1^3 x_2^2+29x_1^2 x_2^3}{5760},
\end{equation*}
which is also correct.
\end{example}

Let $T\in {\rm WMB}(g,n)$ and $\vec u=(u_1,\dots,u_n)$ is an
$n$-vector of nonnegative integers with $|u|=u_1+\cdots+u_n=3g-2+n$.
Since ${\rm WMB}(g,1)$ contains only one element, in the following
discussion, we assume $n\geq2$ with obvious modifications for the
case $n=1$.

We define $P_T(\vec u)$ to be the set of maps $p$ from
$V(T)\backslash L(T)$ to the set of $n$-vector of nonnegative
integers with additional requirements $|{\vec p}_{r_T}|=3g(r_T)$,
where $\vec p_v$ denotes the image of $v$ under $p$ and
$$\vec p_{r_T}+\sum_{\substack{v\notin L(T)\\v\neq r_T}}\vec p_{v}=\vec u-(3g(\ell^{-1}(1)),\dots,3g(\ell^{-1}(n))),$$
where $\ell$ is the bijective labeling map from $L(T)$ to
$\{1,\dots,n\}$. An obvious necessary condition for $P_T(\vec u)$ to
be nonempty is that $u_i\geq 3g(\ell^{-1}(i))$ for all $1\leq i\leq
n$.

\begin{corollary}  \label{tree2}
Let $n\geq2$, $d_i\geq0, \sum_{i=1}^n d_i=3g-3+n$. Then
\begin{multline}
\langle\tau_{d_1}\cdots\tau_{d_n}\rangle_g=\frac{1}{12^g}\sum_{m\geq0}(-1)^m\sum_{\substack{(d_2,\cdots,d_n)=\vec s+\vec t\\ |\vec s|=m}} \binom{m}{\vec s}
\sum_{T\in\ {\rm WMB}(g,n)} \frac{(n-3+2g-2g(r_T))!!}{(n-1+2g)!!}\\
\times \prod_{\substack{v\in L(T)\\v\neq r_T}}\frac{1}{(2g(v))!!}
\prod_{\substack{v\notin
L(T)\\v\neq r_T}}\frac{\left(|L(v)|-3+\sum\limits_{\substack{w\in D(v)\\w\neq
v}}2g(w)\right)!!}{\left(|L(v)|-1+\sum\limits_{w\in
D(v)}2g(w)\right)!!}\\
\times\sum_{p\in P_T(m+1+d_1,\vec t)}  \binom{3g(r_T)}{\vec p_{r_T}}\prod_{\substack{v\notin L(T)\\v\neq r_T}}\binom{3g(v)+1}{\vec p_v}.
\end{multline}
\end{corollary}
\begin{proof} From the proof of Theorem \ref{tree}, we have
\begin{multline} \label{comb2}
(x_1+\dots+x_n) F_g(x_1,\dots,x_n)\\
=\frac{1}{12^g}\sum_{T\in\ {\rm WMB}(g,n)}
\frac{(n-3+2g-2g(r_T))!!}{(n-1+2g)!!}(x_1+\cdots+x_n)^{3g(r_T)}\prod_{\substack{v\in
L(T)\\v\neq r_T}}\frac{x_{\ell(v)}^{3g(v)}}{(2g(v))!!}\\
\times\prod_{\substack{v\notin
L(T)\\v\neq r_T}}\frac{\left(|L(v)|-3+\sum\limits_{\substack{w\in D(v)\\w\neq
v}}2g(w)\right)!!\left(\sum\limits_{w\in
L(v)}x_{\ell(w)}\right)^{3g(v)+1}}{\left(|L(v)|-1+\sum\limits_{w\in
D(v)}2g(w)\right)!!}.
\end{multline}

We may multiply
$$\frac{1}{x_1+\cdots+x_n}=\sum_{m\geq0}(-1)^m\frac{(x_2+\cdots+x_n)^{m}}{x_1^{m+1}}$$
to the right hand side of equation \eqref{comb2}. Comparing coefficients of both sides gives the desired identity.
\end{proof}

\begin{example}
Let us compute $\langle\tau_3\tau_2\rangle_2$ using Corollary
\ref{tree2}. The common factor $1/12^g$ will be counted at last.

\begin{enumerate}
\item[(1)]
When $m=0,\ t=2$. The following two trees have
$P_T(4,2)\neq\emptyset$

\begin{minipage}{0.9in}
\begin{displaymath}
  \vcenter{\tiny\xymatrix@C=1mm@R=0mm{
   &2& \\
    &  *={\bullet}\ar@{-}[ddl]\ar@{-}[ddr] &  \\
  & & \\
 *=<12pt>[o][F-]{1}  &   &
  *=<12pt>[o][F-]{2}  \\
   0&&0
         }}\,
\end{displaymath}
$$\vec p_{r_T}=(4,2)$$
\end{minipage}
\begin{minipage}{0.9in}
\begin{displaymath}
  \vcenter{\tiny\xymatrix@C=1mm@R=0mm{
   &1& \\
    &  *={\bullet}\ar@{-}[ddl]\ar@{-}[ddr] &  \\
  & & \\
 *=<12pt>[o][F-]{1}  &   &
  *=<12pt>[o][F-]{2}  \\
   1&&0
         }}\,
\end{displaymath}
$$\vec p_{r_T}=(1,2)$$
\end{minipage}
\begin{minipage}{2.0in}
$$I_1=\frac{1}{5!!}\binom{6}{4}+\frac{1}{5!!\cdot 2}\binom{3}{1}=\frac{11}{10}.$$
\end{minipage}
\medskip

\item[(2)]
When $m=1,\ t=1$. The following two trees have
$P_T(5,1)\neq\emptyset$

\begin{minipage}{0.9in}
\begin{displaymath}
  \vcenter{\tiny\xymatrix@C=1mm@R=0mm{
   &2& \\
    &  *={\bullet}\ar@{-}[ddl]\ar@{-}[ddr] &  \\
  & & \\
 *=<12pt>[o][F-]{1}  &   &
  *=<12pt>[o][F-]{2}  \\
   0&&0
         }}\,
\end{displaymath}
$$\vec p_{r_T}=(5,1)$$
\end{minipage}
\begin{minipage}{0.9in}
\begin{displaymath}
  \vcenter{\tiny\xymatrix@C=1mm@R=0mm{
   &1& \\
    &  *={\bullet}\ar@{-}[ddl]\ar@{-}[ddr] &  \\
  & & \\
 *=<12pt>[o][F-]{1}  &   &
  *=<12pt>[o][F-]{2}  \\
   1&&0
         }}\,
\end{displaymath}
$$\vec p_{r_T}=(2,1)$$
\end{minipage}
\begin{minipage}{2.0in}
$$I_2=-\frac{1}{5!!}\binom{6}{5}-\frac{1}{5!!\cdot 2}\binom{3}{2}=-\frac{1}{2}.$$
\end{minipage}

\medskip

\item[(3)]
When $m=2,\ t=0$. The following three trees have
$P_T(6,0)\neq\emptyset$

\begin{minipage}{0.8in}
\begin{displaymath}
  \vcenter{\tiny\xymatrix@C=1mm@R=0mm{
   &2& \\
    &  *={\bullet}\ar@{-}[ddl]\ar@{-}[ddr] &  \\
  & & \\
 *=<12pt>[o][F-]{1}  &   &
  *=<12pt>[o][F-]{2}  \\
   0&&0
         }}\,
\end{displaymath}
$$\vec p_{r_T}=(6,0)$$
\end{minipage}
\begin{minipage}{0.8in}
\begin{displaymath}
  \vcenter{\tiny\xymatrix@C=1mm@R=0mm{
   &1& \\
    &  *={\bullet}\ar@{-}[ddl]\ar@{-}[ddr] &  \\
  & & \\
 *=<12pt>[o][F-]{1}  &   &
  *=<12pt>[o][F-]{2}  \\
   1&&0
         }}\,
\end{displaymath}
$$\vec p_{r_T}=(3,0)$$
\end{minipage}
\begin{minipage}{0.8in}
\begin{displaymath}
  \vcenter{\tiny\xymatrix@C=1mm@R=0mm{
   &0& \\
    &  *={\bullet}\ar@{-}[ddl]\ar@{-}[ddr] &  \\
  & & \\
 *=<12pt>[o][F-]{1}  &   &
  *=<12pt>[o][F-]{2}  \\
   2&&0
         }}\,
\end{displaymath}
$$\vec p_{r_T}=(0,0)$$
\end{minipage}
\begin{minipage}{2.0in}
$$I_3=\frac{1}{5!!}+\frac{1}{5!!\cdot 2}+\frac{3!!}{5!!\cdot 4!!}=\frac{1}{8}.$$
\end{minipage}

\end{enumerate}

\medskip

Summing up, we get the desired result

$$\langle\tau_3\tau_2\rangle_2=\frac{I_1+I_2+I_3}{12^2}=\frac{29}{5760}.$$
\end{example}

The authors of the paper \cite{KMZ} proved an explicit formula of
higher Weil-Petersson volumes of moduli spaces of curves in terms of
integrals of $\psi$ classes. For example, in the case of classical
Weil-Petersson volumes, their formula reads

$$\int_{\overline{\mathcal M}_g} \kappa_1^{3g-3}=\sum_{k=1}^{3g-3}\frac{(-1)^{3g-3-k}}{k!}\sum_{\substack {a_1+\cdots+a_k=3g-3\\a_i>0}}\binom{
3g-3}{a_1,\dots,a_k}\langle \tau_{a_1+1}\cdots\tau_{a_k+1}\rangle_g.$$

So via Kaufmann-Manin-Zagier's formula, Corollary \ref{tree2} also
gives a closed formula of higher Weil-Petersson volumes in terms of
summation over WMB trees .

Recall the famous formula of Kontsevich \cite{Ko} expressing
intersection numbers in terms of summation over ribbon graphs
\begin{equation*}
 \sum_{\sum d_i = 3g-3+n} \langle\tau_{d_1}
\dots \tau_{d_n}\rangle \prod_{i=1}^n
\frac{(2d_i-1)!}{\lambda_i^{2d_i+1}} \, =\sum_{\Gamma\in
G^{3}_{g,n}} \frac{2^{-(4g-4+2n)}}{|{\rm Aut}(\Gamma)|} \prod_{e\in
e(\Gamma)} \frac{2}{{\lambda_{1,e}}+{\lambda_{2,e}}},
\end{equation*}
where the summation is over all trivalent ribbon graphs $\Gamma$ of
genus $g$ with $n$ cells, the product is over all edges $e$ of
$\Gamma$, $\lambda_{1,e}$ and $\lambda_{2,e}$ are the $\lambda_i$'s
corresponding to the two sides of an edge $e$.

While the enumeration of ribbon graphs is very difficult, the
enumeration of binary trees is much easier.

Kazarian and Lando \cite{KaL} derived from ELSV formula \cite{ELSV} a closed formula of $\psi$ class integrals in terms of
Hurwitz numbers. Zvonkine \cite{Zv} has an interesting interpretation of the string
and dilaton equations as operations on graphs with marked vertices.

\vskip 30pt
\section{Other applications of $n$-point functions}

We now prove an effective recursion formula which
explicitly expresses intersection indices in terms of intersection
indices with strictly lower genus.

\begin{proposition} \label{eff}
Let $d_j\geq0$ and $\sum_{j=1}^n d_j=3g+n-3$. Then
\begin{multline*}
(2g+n-1)(2g+n-2)\langle\prod_{j=1}^n\tau_{d_j}\rangle_g\\
=\frac{2d_1+3}{12}\langle\tau_0^4\tau_{d_1+1}\prod_{j=2}^n\tau_{d_j}\rangle_{g-1}-\frac{2g+n-1}{6}\langle\tau_0^3\prod_{j=1}^n\tau_{d_j}\rangle_{g-1}\\
+\sum_{\{2,\dots,n\}=I\coprod
J}(2d_1+3)\langle\tau_{d_1+1}\tau_0^2\prod_{i\in
I}\tau_{d_i}\rangle_{g'}\langle\tau_0^2\prod_{i\in
J}\tau_{d_i}\rangle_{g-g'}\\
-\sum_{\{2,\dots,n\}=I\coprod
J}(2g+n-1)\langle\tau_{d_1}\tau_0\prod_{i\in
I}\tau_{d_i}\rangle_{g'}\langle\tau_0^2\prod_{i\in
J}\tau_{d_i}\rangle_{g-g'}.
\end{multline*}
It's not difficult to see that when indices $d_j\geq1$, all non-zero
intesection indices on the right hands have genera strictly less
than $g$.
\end{proposition}
\begin{proof}

First note that Proposition \ref{npoint} (ii) is precisely
$$(2g+n-1)\langle\tau_0\prod_{j=1}^n\tau_{d_j}\rangle_g=\frac{1}{12}\langle\tau_0^4\prod_{j=1}^n\tau_{d_j}\rangle_{g-1}
+\frac{1}{2}\sum_{\underline{n}=I\coprod
J}\langle\tau_0^2\prod_{i\in
I}\tau_{d_i}\rangle_{g'}\langle\tau_0^2\prod_{i\in
J}\tau_{d_i}\rangle_{g-g'}.$$
Applying this, we can group the first and third terms on the
right hand side of Proposition \ref{eff} and further simplify to the following
recursion relation.
\begin{multline*}
(2g+n-1)\langle\tau_r\prod_{j=1}^n\tau_{d_j}\rangle_g=(2r+3)\langle\tau_0\tau_{r+1}\prod_{j=1}^n\tau_{d_j}\rangle_g\\
-\frac{1}{6}\langle\tau_0^3\tau_r\prod_{j=1}^n\tau_{d_j}\rangle_{g-1}
-\sum_{\underline{n}=I\coprod J}\langle\tau_0\tau_r\prod_{i\in
I}\tau_{d_i}\rangle_{g'}\langle\tau_0^2\prod_{i\in
J}\tau_{d_i}\rangle_{g-g'}.
\end{multline*}
So we need only prove the following equivalent statement of
Proposition \ref{eff}:

\begin{multline} \label{ODE2}
y\sum_{g=0}^\infty(2g+n-1)
F_g(y,x_1,\dots,x_n)=2y\frac{\partial}{\partial y}
\left(\left(\sum_{j=1}^n y+x_j\right) F(y,x_1,\dots,x_n)\right)\\
+\left(\left(y+\sum_{j=1}^n x_j\right)-\frac{y}{6}\left(y+\sum_{j=1}^n x_j\right)^3\right) F(y,x_1,\dots,x_n)\\
-y\sum_{\substack{\underline{n}=I\coprod J\\J\neq\emptyset}}
\left(y+\sum_{i\in I} x_i\right)\left(\sum_{i\in J} x_i\right)^2
 F(y,x_I)F(x_J).
\end{multline}

From Witten's ODE \eqref{ODE} in Section 2, it's not difficult to get the following equation.
\begin{multline} \label{ODE3}
2y\left(y+\sum_{j=1}^n x_j\right)\frac{\partial}{\partial
y}\left(\left(y+\sum_{j=1}^n x_j\right) F(y,x_1,\dots,x_n)\right)\\
=\left(\frac{y}{4}\left(y+\sum_{j=1}^n x_j\right)^4-y\left(y+\sum_{j=1}^n x_j\right)-\left(y+\sum_{j=1}^n x_j\right)^2\right) F(y,x_1,\dots,x_n)\\
+y\sum_{\substack{\underline{n}=I\coprod J\\J\neq\emptyset}}
\left(\left(y+\sum_{i\in I} x_i\right)\left(\sum_{i\in J}
x_i\right)^3+2\left(y+\sum_{i\in I} x_i\right)^2\left(\sum_{i\in J}
x_i\right)^2\right) F(y,x_I)F(x_J).
\end{multline}

Multiply each side of the equation \eqref{ODE2} by
$y+\sum_{j=1}^n x_j$ and substitute the differential part using the
above equation \eqref{ODE3}, we get
\begin{multline*}
y\sum_{g=0}^\infty(2g+n-1)\left(y+\sum_{i=1}^n x_i\right) F_g(y,x_1,\dots,x_n)\\
=\left(\frac{y}{12}\left(y+\sum_{j=1}^n x_j\right)^4-y\left(y+\sum_{j=1}^n x_j\right)\right) F(y,x_1,\dots,x_n)\\
+y\sum_{\substack{\underline{n}=I\coprod J\\J\neq\emptyset}}
\left(y+\sum_{i\in I} x_i\right)^2\left(\sum_{i\in J} x_i\right)^2
 F(y,x_I)F(x_J).
\end{multline*}
Add to each side with the term $$y\left(y+\sum_{j=1}^n x_j\right)
F(y,x_1,\dots,x_n),$$ we get the desired equation \eqref{ODE2}. So we
conclude the proof of Proposition \ref{eff}.
\end{proof}

The recursion formula in Proposition \ref{eff}, together with the string and dilaton
equations, provides an effective algorithm for computing
intersection indices on moduli spaces of curves.

Now we prove two interesting combinatorial identities. As pointed
out to us by Lando, these kind of formulae are usually called Abel
identities and they arise naturally in enumeration of various kinds
of marked trees.
\begin{lemma} \label{abel} Let $n\geq2$.
\begin{enumerate}
\item[i)] Assume that if $I=\emptyset$, then $(\sum_{i\in I}x_i)^{|I|}=1$. We have
$$\sum_{\{2,\dots,n\}=I\coprod J}(x_1+\sum_{i\in I}x_i)^{|I|}(-x_1+\sum_{i\in J}x_i)^{|J|}=\sum_{\{2,\dots,n\}=I\coprod J}(\sum_{i\in I}x_i)^{|I|}(\sum_{i\in J}x_i)^{|J|}$$
\item[ii)] We have
\begin{equation*}
\sum_{\underline{n}=I\coprod J\atop I,J\neq\emptyset}(\sum_{i\in I}x_i)^{|I|-1}(\sum_{i\in J}x_i)^{|J|-1}=2(n-1)(\sum_{j=1}^n x_j)^{n-2}
\end{equation*}
\end{enumerate}
\end{lemma}
\begin{proof}
Let $\prod_{j=1}^n x_j^{d_j}$ be any monomial of
\begin{equation}\label{poly}
\sum_{\{2,\dots,n\}=I\coprod J}(x_1+\sum_{i\in
I}x_i)^{|I|}(-x_1+\sum_{i\in J}x_i)^{|J|}.
\end{equation}
Since $\sum_{j=1}^n d_j=n-1$, so if $d_1>0$, then their must exist some $j>1$ such that $d_j=0$.

The statement (i) means that the polynomial \eqref{poly} does not
contain $x_1$, so we need only prove that after substitute $x_n=0$
in \eqref{poly}, the resulting polynomial does not contain $x_1$.
\begin{align*}
&\sum_{\{2,\dots,n-1\}=I\coprod J}\left((x_1+\sum_{i\in I}x_i)^{|I|+1}(-x_1+\sum_{i\in J}x_i)^{|J|}+(x_1+\sum_{i\in I}x_i)^{|I|}(-x_1+\sum_{i\in J}x_i)^{|J|+1}\right)\\
&=(\sum_{j=2}^{n-1}x_j)\sum_{\{2,\dots,n-1\}=I\coprod J}(x_1+\sum_{i\in I}x_i)^{|I|}(-x_1+\sum_{i\in J}x_i)^{|J|}.
\end{align*}
So (i) follows by induction.

We prove statement (ii) by induction. Regard the LHS and RHS of (ii)
as polynomials in $x_n$ with degree $n-2$, we need to prove (ii)
when substitute $x_n=-x_i$ for $i=1\dots n-1$. It's sufficient to
check the case $x_n=-x_{n-1}$.
\begin{align*}
LHS&=2\sum_{\{1,\dots,n-2\}=I\coprod J}\left((x_{n-1}+\sum_{i\in I}x_i)^{|I|}(-x_{n-1}+\sum_{i\in J}x_i)^{|J|}+(\sum_{i\in I}x_i)^{|I|+1}(\sum_{i\in J}x_i)^{|J|-1}\right)\\
&=2\sum_{\{1,\dots,n-2\}=I\coprod J}\left((\sum_{i\in I}x_i)^{|I|}(\sum_{i\in J}x_i)^{|J|}+(\sum_{i\in I}x_i)^{|I|+1}(\sum_{i\in J}x_i)^{|J|-1}\right)\\
&=4(\sum_{j=1}^{n-2} x_j)^{n-2}+(\sum_{j=1}^{n-2}x_j)^2\sum_{\{1,\dots,n-2\}=I\coprod J\atop I,J\neq\emptyset}(\sum_{i\in I}x_i)^{|I|-1}(\sum_{i\in J}x_i)^{|J|-1}\\
&=2(n-1)(\sum_{j=1}^{n-2} x_j)^{n-2}=RHS.
\end{align*}
Note that if a term has power $|J|-1$, then $J\neq\emptyset$ is assumed.
\end{proof}

As an interesting exercise we give a proof of the following
well-known formula.
\begin{corollary}
Let $n\geq3$, $d_j\geq 0$ and $\sum_{j=1}^n d_j=n-3$. Then
$$\langle\tau_{d_1}\cdots\tau_{d_n}\rangle_0=\binom{n-3}{d_1,\dots, d_n}.$$
\end{corollary}
\begin{proof}
It's equivalent to prove that for $n\geq3$
\begin{align*}
(\sum_{j=1}^n x_j)^{n-3}&=G_0(x_1,\dots,x_n)\\
&=\frac{1}{2(n-1)\sum_{j=1}^n x_j}\sum_{\underline{n}=I\coprod J}(\sum_{i\in I}x_i)^2(\sum_{i\in J}x_i)^2G_0(x_I)G_0(x_J)\\
&=\frac{1}{2(n-1)\sum_{j=1}^n x_j}\sum_{\underline{n}=I\coprod J}(\sum_{i\in I}x_i)^{|I|-1}(\sum_{i\in J}x_i)^{|J|-1}.
\end{align*}
This is just Lemma \ref{abel} (ii).
\end{proof}

Finally in this section, we make a remark about the following DVV formula.
\begin{multline*}
\langle\tau_{k+1}\tau_{d_1}\cdots\tau_{d_n}\rangle_g=\frac{1}{(2k+3)!!}\left[\sum_{j=1}^n
\frac{(2k+2d_j+1)!!}{(2d_j-1)!!}\langle\tau_{d_1}\cdots
\tau_{d_{j}+k}\cdots\tau_{d_n}\rangle_g\right.\\
+\frac{1}{2}\sum_{r+s=k-1}
(2r+1)!!(2s+1)!!\langle\tau_r\tau_s\tau_{d_1}\cdots\tau_{d_n}\rangle_{g-1}\\
\left.+\frac{1}{2}\sum_{r+s=k-1} (2r+1)!!(2s+1)!!
\sum_{\underline{n}=I\coprod J}\langle\tau_r\prod_{i\in
I}\tau_{d_i}\rangle_{g'}\langle\tau_s\prod_{i\in
J}\tau_{d_i}\rangle_{g-g'}\right]
\end{multline*}

In fact,  when adopting the conventions \eqref{con} for $F_0(x)$ and $F_0(x,y)$, DVV formula can be written concisely as following:
$$\sum_{j\in\mathbb Z}\ [j-k+\frac{1}{2}]_0^k \left(\mathcal C\left(y^{k-1-j}z^j, F(y,-z,x_1,\dots,x_n)
+\sum_{\underline{n}=I\coprod J}F(y,x_I)F(-z,x_J)\right)\right)=0,$$
where for any integers $m$ and $k\geq -1$,
$$[m+\frac{1}{2}]_0^k=(m+\frac{1}{2})(m+1+\frac{1}{2})\cdots (m+k+\frac{1}{2})=
\frac{\Gamma(m+k+\frac{3}{2})}{\Gamma(m+\frac{1}{2})}.$$

Gathmann \cite{Ga} noticed this fact for the more general Virasoro contraints for Gromov-Witten invariants.

We know there are several proofs that Witten's KdV conjecture implies DVV
formula \cite{DVV, Goe, KaS, La}. See also \cite{Ge}. However we still pose the following problem, which we are not
able to solve for now.

\begin{problem}
Give a direct proof of the above reformulated DVV using the recursion formula of $n$-point function
in Proposition \ref{nptequiv} (ii).
\end{problem}

\vskip 30pt
\section{Vanishing identities}

In fact, the vanishing identities in Proposition \ref{nptcoeff} (ii)  can be generalized
to universal equations for Gromov-Witten invariants \cite{LiuXu3}.

\begin{conjecture} \label{conj1} Let $X$ be a smooth projective variety. Given a basis $\{\gamma_a\}$ for $H^*(X,\mathbb Q)$,
let $x_i, y_i\in H^*(X)$ and $k\geq 2g-3+r+s$. Then the Gromov-Witten potential function satisfies
$$\sum_{g'=0}^g\sum_{j\in\mathbb
Z}(-1)^j\langle\langle\tau_j(\gamma_a)\prod_{i=1}^r\tau_{p_i}(x_i)\rangle\rangle^X_{g'}
\langle\langle\tau_{k-j}(\gamma^a)\prod_{i=1}^s\tau_{q_i}(y_i)\rangle\rangle^X_{g-g'}=0.$$
Note that $j$ runs over all integers and Gathmann's convention \cite{Ga} is used
$$\langle\tau_{-2}(pt)\rangle_{0,0}^X=1$$
and
$$\langle\tau_{m}(\gamma_1)\tau_{-1-m}(\gamma_2)\rangle_{0,0}^X=(-1)^{\max(m,-1-m)}\int_X\gamma_1\cdot\gamma_2,
\quad m\in\mathbb Z.$$ All other Gromov-Witten invariants
containing a negative power of cotangent line classes are defined to be zero.
\end{conjecture}

\begin{conjecture}
Let $k>g$. Then
\begin{equation*}
\sum_{j=0}^{2k}(-1)^j\langle\langle\tau_{j}(\gamma_a)\tau_{2k-j}(\gamma^a)\rangle\rangle_g^X=0.
\end{equation*}
\end{conjecture}

Recently, X. Liu and Pandharipande \cite{Liu,LP} give a proof of the
above conjectures. Their proof uses virtual localization to get
topological recursion relations (TRR) in the tautological ring of
moduli spaces of curves, which are pulled back, via the forgetful
morphism $\pi: \overline{\mathcal M}_{g,n}(X,d)\rightarrow
\overline{\mathcal M}_{g,n}$, to universal equations for
Gromov-Witten invariants. Note that the $\Psi_i$ class on
$\overline{\mathcal M}_{g,n}(X,d)$ differ from $\pi^*(\psi_i)$ by a
cycle containing generic elements whose domain curves consist of one
genus-$g$ and one genus-$0$ components, with the $i$-th marked point
lying on the genus-$0$ component.

For example, one of topological recursion relations proved by X. Liu and Pandharipande \cite{LP}  is the following:
\begin{proposition} (X. Liu-Pandharipande) For $k\geq 2g-1$, there is the following topological relation in $A^{k+1}(\overline{\mathcal M}_{g,2})$
\begin{equation} \label{top}
-\psi_1^{k+1}+(-1)^{k+1}\psi_2^{k+1}+\sum_{\substack{g_1+g_2=g\\i+j=k}}(-1)^i \iota_*(\psi_{*1}^i\psi_{*2}^j\cap[\Delta_{1,2}(g_1,g_2)])=0,
\end{equation}
where $\iota: \Delta_{1,2}\rightarrow \overline{\mathcal M}_{g,2}$ denotes the boundary divisor parametrizing reducible curves $C=C_1\cup C_2$
with markings $p_1\in C_1, p_2\in C_2$ and $C_1\cap C_2=p_*$,  the cotangent line classes of $p_*$ along $C_1$ and $C_2$ are denoted by
$\psi_{*1}, \psi_{*2}$ respectively.

\end{proposition}

The above topological recursion relation \eqref{top} corresponds to the case $r=s=1$ of Conjecture \ref{conj1}, namely for $k\geq 2g-1$
\begin{multline*}
(-1)^{k+1}\langle\langle\tau_{k+q+1}(y)\tau_{p}(x)\rangle\rangle_{g}
-\langle\langle\tau_{k+p+1}(x)\tau_{q}(y)\rangle\rangle_{g}\\
+\sum_{g'=0}^g\sum_{j=0}^{k}(-1)^j\langle\langle\tau_{j}(\gamma_a)\tau_{p}(x)\rangle\rangle_{g'}\langle\langle\tau_{k-j}(\gamma^a)\tau_{q}(y)\rangle\rangle_{g-g'}=0.
\end{multline*}

Similarly TRR also leads to
universal relations for Hodge integrals and Witten's $r$-spin
intersection numbers \cite{Wi2}.

Let $\Sigma$ be a Riemann surface of genus $g$ with marked points
$x_1,x_2,\dots,x_s$. Fix an integer $r\geq 2$. Label each marked
point $x_i$ by an integer $m_i$, $0\leq m_i\leq r-1$. Consider the
line bundle $\mathcal S=K\otimes(\otimes_{i=1}^s\mathcal
O(x_i)^{-m_i})$ over $\Sigma$, where $K$ denotes the canonical line
bundle. If $2g-2-\sum_{i=1}^s m_i$ is divisible by $r$, then there
are $r^{2g}$ isomorphism classes of line bundles $\mathcal T$ such
that $\mathcal T^{\otimes r}\cong \mathcal S$. The choice of an
isomorphism class of $\mathcal T$ determines a moduli space $\mathcal
M^{1/r}_{g,s}$ with compactification $\overline{\mathcal
M}^{1/r}_{g,s}$.

Let $\mathcal V$ be a vector bundle over $\overline{\mathcal
M}^{1/r}_{g,s}$ whose fiber is the dual space to
$H^1(\Sigma,\mathcal T)$. The top Chern class $c_D(\mathcal{V})$ of
this bundle has degree $ D=(g-1)(r-2)/r+\sum_{i=1}^sm_i/r. $

We associate with each marked point $x_i$ an integer $n_i\geq 0$.
Witten's $r$-spin intersection numbers are defined by
$$
\langle\tau_{n_1,m_1}\dots\tau_{n_s,m_s}\rangle_g=
\frac{1}{r^g}\int_{\overline{\mathcal
M}^{1/r}_{g,s}}\prod_{i=1}^s\psi(x_i)^{n_i}\cdot
 c_D(\mathcal{V}),
$$
which is non-zero only if $ (r+1)(2g-2)+rs=r\sum_{j=1}^s
n_j+\sum_{j=1}^s m_j. $

Consider the formal series $F$ in variables $t_{n,m}$, $n\geq 0$ and
$0\leq m\leq r-1$,
$$
F(t_{0,0},t_{0,1},\dots)=\sum_{d_{n,m}}\langle\prod_{n,m}\tau_{n,m}^{d_{n,m}}\rangle
\prod_{n,m}\frac{t_{n,m}^{d_{n,m}}}{d_{n,m}!}.
$$

Let $\eta^{ij}=\delta_{i+j,r-2}$ and
$$\langle\langle\tau_{n_1,m_1}\dots\tau_{n_s,m_s}\rangle\rangle=\frac{\partial}{\partial t_{n_1,m_1}}\dots\frac{\partial}{\partial t_{n_s,m_s}}F(t_{0,0},t_{0,1},\dots).$$

\begin{proposition}
Let $k\geq 2g-3+u+v$ and $0\leq \ell_i,m_i\leq r-2$. Then
$$\sum_{g'=0}^g\sum_{j\in\mathbb
Z}(-1)^j\langle\langle\tau_{j,m'}\prod_{i=1}^u\tau_{p_i,m_i}\rangle\rangle_{g'}
\eta^{m'm''}\langle\langle\tau_{k-j,m''}\prod_{i=1}^v\tau_{q_i,
\ell_i}\rangle\rangle_{g-g'}=0.$$
Note that $j$ runs over all
integers and we define
$$\langle\tau_{-2,r-2}\rangle_0=1$$
and for $0\leq m\leq r-2$,
$$\langle\tau_{n,m}\tau_{-1-n,r-2-m}\rangle_{0}=(-1)^{\max(n,-1-n)},\qquad n\in\mathbb Z.$$
\end{proposition}

\begin{proposition}
Let $k>g$. Then
$$\sum_{j=0}^{2k}(-1)^j\eta^{m'm''}\langle\langle\tau_{j,m'}\tau_{2k-j,m''}\rangle\rangle_g=0.$$
\end{proposition}

These results were also conjectured by us before and follows from X. Liu and Pandharipande's topological recursion relations \cite{LP}.

\vskip 30pt
\appendix
\section{Verification of Witten's ODE}

We will verify that the functions $F_g(y,x_1,\dots,x_n)$ recursively
defined in Proposition \ref{nptequiv} (ii) satisfy Witten's ODE. The proof goes by
inducting on $g$ and $n$, namely we assume $F_h(y,x_1,\dots,x_k)$
satisfies Witten's ODE if either $h<g$ or $k<n$.

Let LHS and RHS denote the left-hand side and right-hand side of the
Witten's ODE. We have
\begin{multline*}
(2g+n)LHS=\frac{y\left(y+\sum_{j=1}^{n}x_j\right)^4}{4}F_{g-1}(y,x_1,\dots,x_n)\\
+\frac{\left(y+\sum_{j=1}^{n}x_j\right)^3}{12}\frac{\partial
}{\partial
y}\left(\left(y+\sum_{j=1}^{n}x_j\right)^2F_{g-1}(y,x_1,\dots,x_n)\right)\\
+y\sum_{\underline{n}=I\coprod J}\left(y+\sum_{i\in
I}x_i\right)^2\left(\sum_{i\in
J}x_i\right)^{2}F_h(y,x_I)F_{g-h}(x_J)\\+\left(y+\sum_{j=1}^{n}x_j\right)y\sum_{\underline{n}=I\coprod
J}\frac{\partial}{\partial y}\left(\left(y+\sum_{i\in
I}x_i\right)^{2}F_h(y,x_i)\right)\left(\sum_{i\in
J}x_i\right)^{2}F_{g-h}(x_J)
\end{multline*}

By induction, we substitute the differential terms using Witten's
ODE and get
\begin{multline*}
(2g+n)LHS=\frac{y\left(y+\sum_{j=1}^{n}x_j\right)^4}{4}F_{g-1}(y,x_{\underline{n}})\\
+\frac{\left(y+\sum_{j=1}^{n}x_j\right)^3}{12}\left(\frac{y}{8}
\left(y+\sum_{j=1}^{n}x_j\right)^4 F_{g-2}(y,x_{\underline{n}})
+\frac{y}{2}\left(y+\sum_{j=1}^{n}x_j\right)F_{g-1}(y,x_1,\dots,x_n)\right.\\
\left.+\frac{y}{2}\sum_{\underline{n}=I\coprod
J}\left(\left(y+\sum_{i\in I}x_i\right)\left(\sum_{i\in
J}x_i\right)^3+2\left(y+\sum_{i\in I}x_i\right)^2\left(\sum_{i\in
J}x_i\right)^2\right)F_h(y,x_I)F_{g-1-h}(x_J)\right.\\
\left.-\frac{1}{2}\left(y+\sum_{j=1}^{n}x_j\right)^{2}F_{g-1}(y,x_1,\dots,x_n)\right)
+y\sum_{\underline{n}=I\coprod J}\left(y+\sum_{i\in
I}x_i\right)^2\left(\sum_{i\in
J}x_i\right)^{2}F_h(y,x_I)F_{g-h}(x_j)\\+\left(y+\sum_{j=1}^{n}x_i\right)\sum_{\underline{n}=I\coprod
J}\left(\frac{y}{8}\left(y+\sum_{i\in
I}x_i\right)^{4}F_{h-1}(y,x_I)+\frac{y}{2}\left(y+\sum_{i\in
I}x_i\right)F_h(y,x_I)\right.\\
\left.+\frac{y}{2}\sum_{{I}=I'\coprod I''}\left(\left(y+\sum_{i\in
I'}x_i\right)\left(\sum_{i\in I''}x_i\right)^3+2\left(y+\sum_{i\in
I'}x_i\right)^2\left(\sum_{i\in
I''}x_i\right)^2\right)F(y,x_I')F(x_I'')\right.\\
\left.-\frac{1}{2}\left(y+\sum_{i\in
I}x_i\right)^2F_h(y,x_I)\right)F_{g-h}(x_J)\left(\sum_{i\in
J}x_i\right)^2
\end{multline*}

Let's introduce some symbols to simplify notations
\begin{align*}
A_g^{a,b}&=\sum_{h=0}^g\sum_{\underline{n}=I\coprod
J}\left(y+\sum_{i\in I}x_i\right)^a\left(\sum_{i\in J}x_i\right)^b
F_h(y,x_I)F_{g-h}(x_J),\\
B_g^{a,b,c}&=\sum_{h=0}^g\sum_{\underline{n}=I\coprod J\coprod
K}\left(y+\sum_{i\in I}x_i\right)^a\left(\sum_{i\in
J}x_i\right)^b\left(\sum_{i\in K}x_i\right)^c
F_h(y,x_I)F_{g-h}(x_J).
\end{align*}
Note that $B_g^{a,b,c}=B_g^{a,c,b}$.

After carefully collecting terms, we arrive at
\begin{multline*}
(2g+n)LHS=\left(\frac{y\left(y+\sum_{j=1}^n
x_j\right)^4}{4}-\frac{\left(y+\sum_{j=1}^n
x_j\right)^4\left(\sum_{j=1}^n
x_j\right)}{24}\right)F_{g-1}(y,x_1,\dots,x_n)\\
+\frac{y\left(y+\sum_{j=1}^n
x_j\right)^7}{96}F_{g-2}(y,x_1,\dots,x_n)+\left(y-\frac{\sum_{j=1}^n
x_j}{2}\right)A^{2,2}_g+\frac{y}{2}A^{1,3}_g \tag{*}\\
+\frac{y}{24}A^{1,6}_{g-1}+\frac{5y}{24}A^{2,5}_{g-1}+\frac{3y}{8}A^{3,4}_{g-1}+\frac{5y}{12}A^{4,3}_{g-1}+\frac{5y}{24}A^{5,2}_{g-1}\\
+\frac{y}{2}B_g^{1,2,4}+\frac{y}{2}B_g^{1,3,3}+\frac{5y}{2}B_g^{2,2,3}+yB_g^{3,2,2}.
\end{multline*}

Substitute the recursion formula for $F_g(x_1,\dots,x_n)$ to the
right-hand side. We have
\begin{multline*}
(2g+n)RHS=\frac{y}{8}\left(y+\sum_{j=1}^{n}x_j\right)^4\left(\frac{\left(y+\sum_{j=1}^{n}x_j\right)^3}{12}F_{g-2}(y,x_1,\dots,x_n)\right.\\
\left.+ \frac{1}{y+\sum_{j=1}^{n}x_j}\sum_{\underline{n}=I\coprod
J}\left(y+\sum_{i\in J}x_i\right)^2\left(\sum_{i\in
J}x_i\right)^{2}F_h(y,x_I)F_{g-1-h}(x_J)\right)\\
+\frac{y}{4}\left(y+\sum_{j=1}^{n}x_j\right)^{4}F_{g-1}(y,x_1,\dots,x_n)+\frac{y}{2}\left(y+\sum_{j=1}^{n}x_j\right)
\left(\frac{\left(y+\sum_{j=1}^{n}x_j\right)^3}{12}F_{g-1}(y,x_1,\dots,x_n)\right.\\
\left.+\frac{1}{y+\sum_{j+1}^{n}x_j}\sum_{h=0}^g\sum_{\underline{n}=I\coprod
J}\left(y+\sum_{i\in I}x_i\right)^2\left(\sum_{i\in
J}x_i\right)^2F_h(y,x_I)F_{g-h}(x_J)\right)\\
+\frac{y}{2}\sum_{h=0}^g\sum_{\underline{n}=I\coprod
J}\left(\left(y+\sum_{i\in I}x_i\right)\left(\sum_{i\in
J}x_i\right)^{3}+2\left(y+\sum_{i\in I}x_i\right)^2\left(\sum_{i\in
J}x_i\right)^2\right)\\
\times \sum_{h=0}^g\left(2h+|I|\right)F_h(y,x_I)F_{g-h}(x_J)\quad \text{(apply Proposition \ref{nptequiv} (ii) to expand)}\\
+\frac{y}{2}\sum_{h=0}^g\sum_{\underline{n}=I\coprod
J}\left(\left(y+\sum_{i\in I}x_i\right)\left(\sum_{i\in
J}x_i\right)^3+2\left(y+\sum_{i\in I}x_i\right)^2\left(\sum_{i\in
J}x_i\right)^2\right)\\
\times \sum_{h=0}^g F_h(y,x_I)(2g-2h+|J|-1)F_{g-h}(x_J)\quad \text{(apply Proposition \ref{nptequiv} (ii) to expand)}\\
+\frac{y}{2}\sum_{h=0}^g\sum_{\underline{n}=I\coprod
J}\left(\left(y+\sum_{i\in I}x_i\right)\left(\sum_{i\in
J}x_i\right)^3+2\left(y+\sum_{i\in I}x_i\right)^2\left(\sum_{i\in
J}x_i\right)^2\right)F_h(y,x_I)F_{g-h}(x_J)\\
-\frac{1}{2}\left(y+\sum_{j=1}^{n}x_j\right)^2\left(\frac{\left(y+\sum_{j=1}^{n}x_j\right)^3}{12}F_{g-1}(y,x_1,\dots,x_n)+\right.\\
\left.\frac{1}{y+\sum_{j=1}^{n}x_j}\sum_{h=0}^g\sum_{\underline{n}=I\coprod
J}\left(y+\sum_{i\in I}x_i\right)^2\left(\sum_{i\in
J}x_i\right)^{2}F_h(y,x_I)F_{g-h}(x_J)\right)
\end{multline*}

After carefully collecting terms, we exactly arrive at
\begin{equation*}
(2g+n)RHS=\text{right hand side of } (*).
\end{equation*}

So we have verified $LHS=RHS$.

$$ \ \ \ \ $$

\end{document}